(Nov. 20, 1995)
# Elementary derivations of summation and transformation formulas for $q$-series

**George Gasper**
Department of Mathematics
Northwestern University
Evanston, IL 60208-2730

We present some elementary derivations of summation and transformation formulas for $q$-series, which are different from, and in several cases simpler or shorter than, those presented in the Gasper and Rahman [1990] "Basic Hypergeometric Series" book (which we will refer to as BHS), the Bailey [1935] and Slater [1966] books, and in some papers; thus providing deeper insights into the theory of $q$-series. Our main emphasis is on methods that can be used to **derive** formulas, rather than to just *verify* previously derived or conjectured formulas. In §5 this approach leads to the derivation of a new family of summation formulas for very-well-poised basic hypergeometric series ${}_{6+2k}W_{5+2k}$, $k = 1, 2, \ldots$. Several of the observations in this paper were presented, along with related exercises, in the author's minicourse on "$q$-Series" at the Fields Institute miniprogram on "Special Functions, $q$-Series and Related Topics," June 12–14, 1995. As is customary, we employ the notations used in BHS for the *shifted factorial*

$$(a)_0 = 1, \ (a)_k = a(a+1)\cdots(a+k-1), \qquad k = 1, 2, \ldots,$$

the *$q$-shifted factorial*

$$(a;q)_0 = 1, \ (a;q)_k = (1-a)(1-aq)\cdots(1-aq^{k-1}), \qquad k = 1, 2, \ldots,$$

$$(a;q)_\infty = \lim_{k\to\infty}(a;q)_k = \prod_{k=0}^{\infty}(1-aq^k), \qquad |q| < 1,$$

$$(a;q)_\alpha = \frac{(a;q)_\infty}{(aq^\alpha;q)_\infty}, \qquad 0 < |q| < 1,$$

the ${}_rF_s$ *hypergeometric series*

$$ {}_rF_s(a_1, a_2, \ldots, a_r; b_1, \ldots, b_s; z) \equiv {}_rF_s\left[\begin{array}{c} a_1, a_2, \ldots, a_r \\ b_1, \ldots, b_s \end{array}; z\right]$$

$$= \sum_{k=0}^{\infty} \frac{(a_1)_k (a_2)_k \cdots (a_r)_k}{k!\, (b_1)_k \cdots (b_s)_k} z^k,$$

---

1991 *Mathematics Subject Classification*. Primary 33D15, 33D20, 33D65; Secondary 33C05, 33C20.

This work was supported in part by the National Science Foundation under grant DMS-9401452.

and the $_r\phi_s$ *basic hypergeometric series*

$$_r\phi_s(a_1, a_2, \ldots, a_r; b_1, \ldots, b_s; q, z) \equiv {_r\phi_s} \begin{bmatrix} a_1, a_2, \ldots, a_r \\ b_1, \ldots, b_s \end{bmatrix}; q, z$$

$$= \sum_{k=0}^{\infty} \frac{(a_1, a_2, \ldots, a_r; q)_k}{(q, b_1, \ldots, b_s; q)_k} \left[(-1)^k q^{\binom{k}{2}}\right]^{1+s-r} z^k,$$

where $\binom{k}{2} = k(k-1)/2$, $(a_1, a_2, \ldots, a_r; q)_k = (a_1; q)_k (a_2; q)_k \cdots (a_r; q)_k$ and the principal value of $q^\alpha$ is taken. We also employ the compact notation

$$_{r+1}W_r(a_1; a_4, a_5, \ldots, a_{r+1}; q, z)$$

for the *very-well-poised* $_{r+1}\phi_r$ series

$$_{r+1}\phi_r \begin{bmatrix} a_1, qa_1^{\frac{1}{2}}, -qa_1^{\frac{1}{2}}, a_4, \ldots, a_{r+1} \\ a_1^{\frac{1}{2}}, -a_1^{\frac{1}{2}}, qa_1/a_4, \ldots, qa_1/a_{r+1} \end{bmatrix}; q, z$$

and define the *bilateral basic hypergeometric* $_r\psi_s$ series by

$$_r\psi_s(z) \equiv {_r\psi_s} \begin{bmatrix} a_1, a_2, \ldots, a_r \\ b_1, b_2, \ldots, b_s \end{bmatrix}; q, z$$

$$= \sum_{k=-\infty}^{\infty} \frac{(a_1, a_2, \ldots, a_r; q)_k}{(b_1, b_2, \ldots, b_s; q)_k} (-1)^{(s-r)k} q^{(s-r)\binom{k}{2}} z^k.$$

For simplicity, unless stated otherwise we shall assume that $n$ is a nonnegative integer, $|q| < 1$ in nonterminating $q$-series, and that the parameters and variables are complex numbers such that the series converge absolutely and any singularities are avoided (which usually leads to isolated conditions on the parameters and variables since the singularities are usually at poles and at limits of sequences of poles). For a discussion of when the above series converge, see Sections 1.2 and 5.1 in BHS (in the third paragraph on p.5 each of the ratios $|b_1 b_2 \cdots b_s|/|a_1 a_2 \cdots a_r|$ should be replaced by $|b_1 b_2 \cdots b_s q|/|a_1 a_2 \cdots a_r|$).

## 1. The $q$-binomial theorem

The summation formula

$$_1F_0(a; -; z) = \sum_{k=0}^{\infty} \frac{(a)_k}{k!} z^k = (1-z)^{-a}, \qquad |z| < 1, \tag{1.1}$$

is called the *binomial theorem* because, when $-a = n$ is a nonnegative integer and $z = -x/y$, it reduces to the binomial theorem for the $n$-th power of the binomial $x + y$:

$$(x+y)^n = \sum_{k=0}^{n} \binom{n}{k} x^k y^{n-k}. \tag{1.2}$$

Since, by l'Hôpital's rule,

$$\lim_{q \to 1} \frac{1-q^a}{1-q} = a$$

and hence

$$\lim_{q \to 1} \frac{(q^a; q)_k}{(q; q)_k} = \frac{(a)_k}{k!},$$

it is natural to consider what happens when the coefficient $(a)_k/k!$ of $z^k$ in the series in (1.1) is replaced by $(q^a;q)_k/(q;q)_k$ or, more generally, by $(a;q)_k/(q;q)_k$. Hence, let us set

$$f(a,z) = \sum_{k=0}^{\infty} \frac{(a;q)_k}{(q;q)_k} z^k, \qquad |z| < 1, \tag{1.3}$$

with $|q| < 1$. The case when $|q| > 1$ will be considered later. Note that, by the Weierstrass M-test, since $|q| < 1$ the series in (1.3) converges uniformly on compact subsets of the unit disk $\{z : |z| < 1\}$ to a function $f(a,z)$ that is an analytic function of $z$ (and of $a$) when $|z| < 1$. One way to find a formula for $f(a,z)$ that is a generalization of (1.1) is to first observe that, since $1 - a = 1 - aq^k + aq^k - a = (1 - aq^k) - a(1 - q^k)$,

$$\begin{aligned}
f(a,z) &= 1 + \sum_{k=1}^{\infty} \frac{(a;q)_k}{(q;q)_k} z^k \\
&= 1 + \sum_{k=1}^{\infty} \frac{(aq;q)_{k-1}}{(q;q)_k} [(1 - aq^k) - a(1 - q^k)] z^k \\
&= 1 + \sum_{k=1}^{\infty} \frac{(aq;q)_k}{(q;q)_k} z^k - a \sum_{k=1}^{\infty} \frac{(aq;q)_{k-1}}{(q;q)_{k-1}} z^k \\
&= f(aq,z) - azf(aq,z) = (1 - az)f(aq,z).
\end{aligned} \tag{1.4}$$

By iterating this functional equation $n - 1$ times, we find that

$$f(a,z) = (az;q)_n \ f(aq^n, z),$$

which on letting $n \to \infty$ and using $q^n \to 0$ yields

$$f(a,z) = (az;q)_\infty \ f(0,z). \tag{1.5}$$

Now set $a = q$ in (1.5) to get

$$f(0,z) = \frac{f(q,z)}{(qz;q)_\infty} = \frac{(1-z)^{-1}}{(qz;q)_\infty} = \frac{1}{(z;q)_\infty},$$

which, combined with (1.3) and (1.5), gives the *q-binomial theorem*

$${}_1\phi_0(a;-;q,z) = \sum_{k=0}^{\infty} \frac{(a;q)_k}{(q;q)_k} z^k = \frac{(az;q)_\infty}{(z;q)_\infty}, \qquad |z| < 1, \ |q| < 1. \tag{1.6}$$

This summation formula was derived by Cauchy [1843], Jacobi [1846], and Heine [1847]. Heine's proof of (1.6), which is reproduced in the books Heine [1878], Bailey [1935, p. 66], Slater [1966, p. 92], and in §1.3 of BHS along with some motivation from Askey [1980], consists of using series manipulations to derive the functional equation

$$(1-z)f(a,z) = (1-az)f(a,qz), \tag{1.7}$$

iterating (1.7) $n - 1$ times, and then letting $n \to \infty$ to get

$$f(a,z) = \frac{(az;q)_n}{(z;q)_n} f(a, q^n z) = \frac{(az;q)_\infty}{(z;q)_\infty} f(a,0) = \frac{(az;q)_\infty}{(z;q)_\infty},$$

which gives (1.6).

Another derivation of the $q$-binomial theorem can be given by calculating the coefficients $c_k = g_a^{(k)}(0)/k!$, $k = 0, 1, 2, \ldots$, in the Taylor series expansion of the function

$$g_a(z) = \frac{(az;q)_\infty}{(z;q)_\infty} = \sum_{k=0}^\infty c_k z^k, \tag{1.8}$$

which is an analytic function of $z$ when $|z| < 1$ and $|q| < 1$. Clearly $c_0 = g_a(0) = 1$. One may show that $c_1 = g_a'(0) = (1-a)/(1-q)$ by taking the logarithmic derivative of $(az;q)_\infty/(z;q)_\infty$ and then setting $z = 0$. But, unfortunately, the succeeding higher order derivatives of $g_a(z)$ become more and more difficult to calculate for $|z| < 1$, and so one is forced to abandon this approach and to search for another way to calculate all of the $c_k$ coefficients. One simple method is to notice that from the definition of $g_a(z)$ as the quotient of two infinite products, it immediately follows that $g_a(z)$ satisfies the functional equation

$$(1 - z)\, g_a(z) = (1 - az)\, g_a(qz), \tag{1.9}$$

which is of course the same as the functional equation (1.7) satisfied by $f(a,z)$. In a *verification* type proof of the $q$-binomial theorem, (1.9) provides *substantial* motivation for showing, as in Heine's proof, that the sum of the $q$-binomial series $f(a,z)$ satisfies the functional equation (1.7).

To calculate the $c_k$ coefficients, we first use (1.9) to obtain

$$\sum_{k=0}^\infty c_k\, z^k - \sum_{k=0}^\infty c_k\, z^{k+1} = \sum_{k=0}^\infty c_k\, q^k\, z^k - a \sum_{k=0}^\infty c_k\, q^k\, z^{k+1},$$

or, equivalently,

$$1 + \sum_{k=1}^\infty (c_k - c_{k-1})\, z^k = 1 + \sum_{k=1}^\infty (c_k q^k - a c_{k-1}\, q^{k-1})\, z^k,$$

which implies that

$$c_k - c_{k-1} = c_k q^k - a c_{k-1}\, q^{k-1}$$

and hence

$$c_k = \frac{1 - aq^{k-1}}{1 - q^k}\, c_{k-1}, \qquad k = 1, 2, \ldots. \tag{1.10}$$

Iterating the recurrence relation (1.10) gives

$$c_k = \frac{(a;q)_k}{(q;q)_k} c_0 = \frac{(a;q)_k}{(q;q)_k}, \qquad k = 0, 1, 2, \ldots,$$

and concludes the derivation of (1.6). For a combinatorial proof of the $q$-binomial theorem using a bijection between two classes of partitions, see Andrews [1969].

It is of interest to note that if $|q| > 1$, then by replacing $q$ in (1.6) by $q^{-1}$, applying the *inversion identity*

$$(a;q)_k = (a^{-1};q^{-1})_k (-a)^k q^{\binom{k}{2}}, \qquad k = 0, 1, 2, \ldots, \tag{1.11}$$

replacing $a$ by $a^{-1}$ and then $z$ by $az/q$, it follows that when $|q| > 1$ the $q$-binomial theorem takes the form

$${}_1\phi_0(a; \text{---}; q, z) = \frac{(z/q; q^{-1})_\infty}{(az/q; q^{-1})_\infty}, \qquad |az/q| < 1,\ |q| > 1. \tag{1.12}$$

In the special case when $z = -x/y$ and $a = q^{-n}$, $n = 0, 1, 2, \ldots$, both (1.6) and (1.12) give that (1.2) has a $q$-analogue of the form

$$y^n(-xq^{-n}/y;q)_n = \sum_{k=0}^{n} \frac{(q^{-n};q)_k}{(q;q)_k}(-1)^k x^k y^{n-k}. \tag{1.13}$$

Another $q$-analogue of (1.2), which we will utilize in the next section, may be easily derived by first observing that from (1.6) we have the product formula

$$\sum_{j=0}^{\infty} \frac{(a;q)_j}{(q;q)_j} z^j \sum_{k=0}^{\infty} \frac{(b;q)_k}{(q;q)_k}(az)^k = \sum_{n=0}^{\infty} \frac{(ab;q)_n}{(q;q)_n} z^n, \qquad |z| < 1, \tag{1.14}$$

which is a $q$-analogue of $(1-z)^{-a}(1-z)^{-b} = (1-z)^{-a-b}$. Then set $j = n - k$ in the product on the left side of (1.14) and compare the coefficients of $z^n$ on both sides of the equation to get

$$\frac{(ab;q)_n}{(q;q)_n} = \sum_{k=0}^{n} \frac{(a;q)_{n-k}(b;q)_k}{(q;q)_{n-k}(q;q)_k} a^k, \tag{1.15}$$

which gives a $q$-analogue of (1.2) in the form

$$(ab;q)_n = \sum_{k=0}^{n} \begin{bmatrix} n \\ k \end{bmatrix}_q (a;q)_{n-k}(b;q)_k a^k, \tag{1.16}$$

where the *q-binomial coefficient* is defined by

$$\begin{bmatrix} n \\ k \end{bmatrix}_q = \frac{(q;q)_n}{(q;q)_k(q;q)_{n-k}}, \qquad k = 0, 1, \ldots, n.$$

Replacing $a$ in (1.16) by $q^{1-n}/c$ and manipulating the $q$-shifted factorials via the identities (I.8), (I.14), and (I.42) in Appendix I of BHS shows that (1.16) is equivalent to the *q-Chu–Vandermonde summation formula*

$$_2\phi_1(q^{-n}, b; c; q, q) = \frac{(c/b;q)_n}{(c;q)_n} b^n, \tag{1.17}$$

which is a $q$-analogue of the *Chu–Vandermonde summation formula* (see p. 2 in BHS)

$$F(-n, b; c; 1) = \frac{(c-b)_n}{(c)_n}.$$

By either switching the order of summation or inverting the base $q$ via (1.11), we find that (1.17) is also equivalent to the summation formula

$$_2\phi_1(q^{-n}, b; c; q, cq^n/b) = \frac{(c/b;q)_n}{(c;q)_n}. \tag{1.18}$$

## 2. Analytic continuations

By manipulating the products on the right side of (1.18) we find that

$$_2\phi_1(z^{-1}, b; c; q, cz/b) = \frac{(c/b, cz;q)_\infty}{(c, cz/b;q)_\infty} = \prod_{k=0}^{\infty} \frac{(1-cq^k/b)(1-czq^k)}{(1-cq^k)(1-czq^k/b)} \tag{2.1}$$

with $z = q^n$, $n = 0, 1, \ldots$, and $|q| < 1$. The infinite product on the right side of (2.1) converges uniformly on compact subsets of the disk $\{z : |cz/b| < 1\}$ to an analytic function of $z$. In view of the identity

$$(z^{-1}; q)_k z^k = (z - 1)(z - q) \cdots (z - q^{k-1}), \qquad k = 0, 1, \ldots,$$

the series on the left side of (2.1) is a sum of analytic functions of $z$, which also converges uniformly on compact subsets of the disk $\{z : |cz/b| < 1\}$ to an analytic function of $z$. Since $q^n \to 0$ as $n \to \infty$ when $|q| < 1$, it follows by analytic continuation that (2.1) holds for arbitrary complex values of $z$ when $|cz/b| < 1$ and $|q| < 1$. Hence, by setting $z = a^{-1}$ we have derived the Jacobi [1846] and Heine [1847] *q-Gauss summation formula*

$$_2\phi_1(a, b; c; q, c/ab) = \frac{(c/a, c/b; q)_\infty}{(c, c/ab; q)_\infty}, \qquad |c/ab| < 1. \tag{2.2}$$

Similarly, from (1.17)

$$_2\phi_1(a, b; c; q, q) = \frac{(a/c, b/c; q^{-1})_\infty}{(1/c, ab/c; q^{-1})_\infty}, \qquad |q| > 1, \tag{2.3}$$

with $a = q^{-n}$, $n = 0, 1, \ldots$, and, because $q^{-n} \to 0$ as $n \to \infty$ when $|q| > 1$, it follows by analytic continuation that (2.3) holds for arbitrary complex values of $a$ when, for convergence, $|ab/c| < 1$. Formula (2.3) also follows from (2.2) by inverting the base $q$.

To see that both (2.2) and (2.3) are $q$-analogues of Gauss' [1813] famous summation formula

$$_2F_1(a, b; c; 1) = \frac{\Gamma(c)\Gamma(c - a - b)}{\Gamma(c - a)\Gamma(c - b)}, \qquad \text{Re}\,(c - a - b) > 0, \tag{2.4}$$

it suffices to replace $a, b, c$ by $q^a, q^b, q^c$, respectively, with $0 < q < 1$ in (2.2) and $q > 1$ in (2.3), and then let $q \to 1$. The $q$-Gauss formula (2.2) was derived in BHS by using the $q$-binomial theorem to derive Heine's [1847] $_2\phi_1$ transformation formula

$$_2\phi_1(a, b; c; q, z) = \frac{(b, az; q)_\infty}{(c, z; q)_\infty} \,_2\phi_1(c/b, z; az; q, b), \qquad |z| < 1, \ |b| < 1, \tag{2.5}$$

setting $z = c/ab$ to reduce the series on the right side of this transformation formula to a $_1\phi_0$ series, and then summing that series via the $q$-binomial theorem.

If we shift the index of summation $k$ in (1.6) by replacing it by $k + n$, we obtain that if $0 < |q| < 1$, $0 < |z| < 1$, and $n = 0, 1, 2, \ldots$, then

$$\frac{(az; q)_\infty}{(z; q)_\infty} = \sum_{k=-\infty}^{\infty} \frac{(a; q)_k}{(q; q)_k} z^k = \sum_{k=-\infty}^{\infty} \frac{(a; q)_{k+n}}{(q; q)_{k+n}} z^{k+n}$$

$$= \sum_{k=-\infty}^{\infty} \frac{(a; q)_n (aq^n; q)_k}{(q; q)_n (q^{n+1}; q)_k} z^{k+n} = \frac{(a; q)_n}{(q; q)_n} z^n \sum_{k=-\infty}^{\infty} \frac{(aq^n; q)_k}{(q^{n+1}; q)_k} z^k, \tag{2.6}$$

where, as usual, the definition of $(a; q)_k$ is extended to negative integer values of $k$ by defining

$$(a; q)_{-k} = \frac{1}{(aq^{-k}; q)_k} = \frac{(-q/a)^k}{(q/a; q)_k} q^{\binom{k}{2}}, \qquad k = 0, 1, 2, \ldots.$$

After replacing $a$ by $aq^{-n}$ and then setting $q^{n+1} = b$, the left and right sides of (2.6) give

$$_1\psi_1(a;b;q,z) = \sum_{k=-\infty}^{\infty} \frac{(a;q)_k}{(b;q)_k} z^k = \frac{(q,b/a,az,q/az;q)_\infty}{(b,q/a,z,b/az;q)_\infty} \qquad (2.7)$$

for $b = q^{n+1}$ when $0 < |q| < 1, 0 < |z| < 1$, and $n = 0, 1, 2, \ldots$. Since $q^{n+1} \to 0$ as $n \to \infty$ when $|q| < 1$, and the infinite series and the infinite product on the left and right sides, respectively, of (2.7) converge to analytic functions of $b$ when $|b| < \min(1, |az|)$ and $|z| < 1$, it follows by analytic continuation that we have derived *Ramanujan's $_1\psi_1$ summation formula*

$$_1\psi_1(a;b;q,z) = \frac{(q,b/a,az,q/az;q)_\infty}{(b,q/a,z,b/az;q)_\infty}, \qquad |b/a| < |z| < 1, \qquad (2.8)$$

which reduces to the $q$-binomial theorem when $b = q$.

The above derivation of (2.8) is essentially in the reverse order of Ismail's [1977] proof of (2.8), which first *reduces* the proof of (2.8) to the case when $b = q^{n+1}$, where $n$ is a nonnegative integer, and then verifies this case by using a shift in the index of summation to obtain a series that is summable by the $q$-binomial theorem. For other proofs of (2.8) and historical comments, see Berndt [1993], BHS, and their references. It should be noted that if we set $b = 0$ in (2.8), replace $q$ and $z$ by $q^2$ and $-qz/a$, respectively, and then let $a \to \infty$, we obtain Jacobi's [1829] *triple product identity*

$$\sum_{k=-\infty}^{\infty} q^{k^2} z^k = \left(q^2, -qz, -q/z; q^2\right)_\infty. \qquad (2.9)$$

The infinite product representations for the *theta functions* $\vartheta_1(x), \vartheta_2(x), \vartheta_3(x)$, and $\vartheta_4(x)$ displayed on p. 13 of BHS are special cases of (2.9).

It is natural to investigate what happens when the shift in index of summation method is applied to the $q$-Gauss summation formula (2.2). Proceeding as in (2.6) with the parameters $a, b, c$ in (2.2) replaced by $aq^{-n}, bq^{-n}, cq^{-n}$, respectively, we find that if $|q| < 1$ and $d = q^{n+1}$, $n = 0, 1, 2, \ldots$, then

$$_2\psi_2(a,b;c,d;q,cd/abq) = \frac{d(q,c/a,c/b,d/a,d/b;q)_\infty}{qc^n(q/a,q/b,c,d,cd/abq;q)_\infty}, \qquad (2.10)$$

when $0 < |cd/abq| < 1$ and $c \neq q^{n-k}$, $k = 0, 1, 2, \ldots$. Unfortunately, the right side of (2.10) is not an analytic function of $d$ in a neighborhood of the origin because of the $c^n$ factor and the initial $d = q^{n+1}$ condition, and the term with index $k$ in the series on the left side of (2.10) has a pole of order $-k$ at $d = 0$ when $k$ is a negative integer. Therefore we cannot analytically continue (2.10) in $d$ to derive an infinite product representation for the series on the left side of (2.10) that is valid for $d$ in a neighborhood of the origin. This helps to explain why *no one* has been able to extend the $q$-Gauss summation formula (2.2) to derive an infinite product representation for a $_2\psi_2$ series that is a $q$-analogue of Dougall's [1907] bilateral hypergeometric series summation formula

$$\sum_{k=-\infty}^{\infty} \frac{(a)_k (b)_k}{(c)_k (d)_k} = \frac{\Gamma(c)\Gamma(d)\Gamma(1-a)\Gamma(1-b)\Gamma(c+d-a-b-1)}{\Gamma(c-a)\Gamma(c-b)\Gamma(d-a)\Gamma(d-b)}, \qquad (2.11)$$

where Re $(c + d - a - b - 1) > 0$ and $(a)_k = (-1)^k/(1-a)_{-k}$, $k = -1, -2, \ldots$, which reduces to Gauss' formula (2.4) when $d = 1$. However, see Ex. 5.20 in BHS for two transformation formulas for $_2\psi_2(a, b; c, d; q, z)$ series.

### 3. The $q$-Pfaff–Saalschütz summation formula

In order to extend the $q$-Chu–Vandermonde and $q$-Gauss summation formulas to $_3\phi_2$ series, we use (1.17) in the form

$$\frac{(a;q)_k}{(c;q)_k} = \frac{(a;q)_n}{(c;q)_n} {}_2\phi_1(q^{k-n}, c/a; q^{1-n}/a; q, q), \qquad 0 \leq k \leq n,$$

to obtain

$$\sum_{k=0}^{n} \frac{(q^{-n}, a, b; q)_k}{(q, c, d; q)_k} q^k = \frac{(a;q)_n}{(c;q)_n} \sum_{k=0}^{n} \frac{(q^{-n}, b; q)_k}{(q, d; q)_k} q^k \sum_{j=0}^{n-k} \frac{(q^{k-n}, c/a; q)_j}{(q, q^{1-n}/a; q)_j} q^j$$

$$= \frac{(a;q)_n}{(c;q)_n} \sum_{j=0}^{n} \frac{(q^{-n}, c/a; q)_j}{(q, q^{1-n}/a; q)_j} {}_2\phi_1(q^{j-n}, b; d; q, q)$$

by a change in order of summation, which on using (1.17) to sum the $_2\phi_1$ series gives the $_3\phi_2$ transformation formula

$$\begin{aligned}&{}_3\phi_2(q^{-n}, a, b; c, d; q, q) \\ &= \frac{(a, d/b; q)_n}{(c, d; q)_n} b^n {}_3\phi_2(q^{-n}, c/a, q^{1-n}/d; q^{1-n}/a, q^{1-n}b/d; q, q).\end{aligned} \qquad (3.1)$$

This formula is a generalization of Heine's $_2\phi_1$ transformation formula (2.5), which follows from (3.1) by setting $d = q^{1-n}/z$, letting $n \to \infty$, and then switching $a$ and $b$. Now notice that if $c/a = q^{1-n}b/d$, i.e. if the *balanced* condition $cd = abq^{1-n}$ holds, then the $_3\phi_2$ series on the right side of (3.1) reduces to a terminating $_2\phi_1$ series that can be summed by (1.17) to give Jackson's [1910] summation formula for a terminating balanced $_3\phi_2$ series

$$_3\phi_2(a, b, q^{-n}; c, abc^{-1}q^{1-n}; q, q) = \frac{(c/a, c/b; q)_n}{(c, c/ab; q)_n}. \qquad (3.2)$$

This formula is usually called the *$q$-Saalschütz* or the *$q$-Pfaff–Saalschütz summation formula* because it is a $q$-analogue of the Pfaff [1797] and Saalschütz [1890] summation formula

$$_3F_2(a, b, -n; c, 1 + a + b - c - n; 1) = \frac{(c-a)_n(c-b)_n}{(c)_n(c-a-b)_n}.$$

Notice that letting $a \to 0$ in (3.2) gives (1.17), letting $a \to \infty$ gives (1.18), and letting $n \to \infty$ in (3.2) gives the $q$-Gauss summation formula (2.2).

If, as in our derivation of the $q$-Gauss formula (2.2), we observe that (3.2) can be rewritten in the form

$$_3\phi_2(a, b, 1/z; c, abq/cz; q, q) = \frac{(c/a, c/b, cz, cz/ab; q)_\infty}{(c, c/ab, cz/a, cz/b; q)_\infty}, \qquad |q| < 1, \qquad (3.3)$$

with $z = q^n$, $n = 0, 1, \ldots$, then the infinite product on the right side of (3.3) clearly converges to an analytic function of $z$ when $|z| < \min(|a/c|, |b/c|)$. However, because $(1/z; q)_k/(abq/cz; q)_k$, $k = 1, 2, \ldots$, has poles at the points $z = abq^j/c$, $1 \leq j \leq k$, there is no neighborhood of the origin in which all of the terms of the series on the left side of (3.3) are analytic; thus (3.3) cannot be analytically continued in $z$ to a neighborhood of the origin. Nevertheless, see equation (2.10.12) in BHS for a nonterminating extension of (3.2) with the sum of *two* balanced $_3\phi_2$ series, and also the bilateral nonterminating extension in Ex. 5.13 of BHS.

The $q$-Pfaff–Saalschütz formula (3.2) was derived in Bailey [1935] and in Slater [1966] by first using an induction argument to give Jackson's [1921] *verification* type proof of his summation formula for a terminating $_8W_7$ series (derived in §5)

$$_8W_7(a;b,c,d,e,q^{-n};q,q) = \frac{(aq,aq/bc,aq/bd,aq/cd;q)_n}{(aq/b,aq/c,aq/d,aq/bcd;q)_n}, \tag{3.4}$$

where $bcde = a^2 q^{n+1}$, and then replacing $d$ by $aq/d$ in (3.4) and letting $a \to \infty$ to get (3.2) as a limit case of (3.4). In BHS, (3.2) was derived by iterating (2.5) twice to derive another of Heine's [1847] transformation formulas

$$_2\phi_1(a,b;c;q,z) = \frac{(abz/c;q)_\infty}{(z;q)_\infty} \, _2\phi_1(c/a,c/b;c;q,abz/c), \tag{3.5}$$

where $|z| < 1$ and $|abz/c| < 1$, using the $q$-binomial theorem to expand the coefficient of the $_2\phi_1$ series on the right side of (3.5) as a power series in powers of $z$, and then equating the coefficients of $z^n$ on both sides of the resulting formula to get (3.2). Recently, Ismail [1995] used the Askey-Wilson difference operators to derive (3.2). For some combinatorial proofs of (3.2), see Andrews and Bressoud [1984], Goulden [1985], and Zeilberger [1987]. Also see Wilf and Zeilberger [1992] for a computer-constructible WZ proof of (3.2).

## 4. Summation formulas for some very-well-poised series

Let us start by deriving some transformation formulas for $_4W_3(a;b;q,z)$ series. Let $|q| < 1$. Since

$$(1-a)(1-b)\frac{(qa^{\frac{1}{2}},-qa^{\frac{1}{2}};q)_k}{(a^{\frac{1}{2}},-a^{\frac{1}{2}};q)_k} = (1-b)(1-aq^{2k}) = 1 + abq^{2k} - (b+aq^{2k})$$
$$= (1 + abq^{2k} - aq^k - bq^k) - (b + aq^{2k} - aq^k - bq^k)$$
$$= (1-aq^k)(1-bq^k) - b(1-aq^k/b)(1-q^k),$$

we have

$$_4W_3(a;b;q,z) = 1 + \sum_{k=1}^\infty \frac{(aq,bq;q)_{k-1}}{(q,aq/b;q)_k}(1-b)(1-aq^{2k})z^k$$
$$= 1 + \sum_{k=1}^\infty \frac{(aq,bq;q)_k}{(q,aq/b;q)_k} z^k - bz \sum_{k=1}^\infty \frac{(aq,bq;q)_{k-1}}{(q,aq/b;q)_{k-1}} z^{k-1}$$
$$= (1-bz) \, _2\phi_1(aq,bq;aq/b;q,z), \quad |z| < 1. \tag{4.1}$$

Hence $_4W_3(a;b;q,b^{-1}) = 0$ when $|b^{-1}| < 1$, and it follows from the special cases $b = q^{-n}$, $n = 1, 2, \ldots$, and the fact that $_4W_3(a;1;q,z) = 0$ that

$$_4W_3(a;q^{-n};q,q^n) = \delta_{n,0}, \tag{4.2}$$

where $\delta_{n,m}$ is the Kronecker delta function. Some recent multidimensional generalizations of (4.2) are given in Bhatnagar and Milne [1995].

This derivation of (4.2) is substantially simpler than that in §2.3 of BHS, which used the $q$-Pfaff–Saalschütz formula (3.2) and the Bailey [1941] and Daum [1942] *q-Kummer summation formula*

$$_2\phi_1(a,b;aq/b;q,-q/b) = \frac{(-q;q)_\infty (aq, aq^2/b^2;q^2)_\infty}{(aq/b, -q/b;q)_\infty}, \qquad |q/b| < 1, \qquad (4.3)$$

derived in §1.8 of BHS. Formula (4.2) can also be derived by the finite difference method employed in Gasper [1989] and Gasper and Rahman [1990a] to derive bibasic extensions of (4.2), and by the method pointed out in Rahman [1990].

Using (4.2), we obtain the expansion formula

$$u_0 = \sum_{k\geq 0} u_k \,\delta_{k,0} = \sum_{k\geq 0} u_k \sum_{j=0}^k \frac{(a, qa^{\frac{1}{2}}, -qa^{\frac{1}{2}}, q^{-k};q)_j}{(q, a^{\frac{1}{2}}, -a^{\frac{1}{2}}, aq^{k+1};q)_j} q^{jk}$$

$$= \sum_{j\geq 0} \frac{(a, qa^{\frac{1}{2}}, -qa^{\frac{1}{2}};q)_j}{(a^{\frac{1}{2}}, -a^{\frac{1}{2}}, aq^{j+1};q)_j} (-1)^j q^{\binom{j}{2}} \sum_{k\geq 0} \frac{(q^{j+1}, aq^{j+1};q)_k}{(q, aq^{2j+1};q)_k} u_{j+k} \qquad (4.4)$$

by a change in order of summation, where $\{u_k\}$ is a sequence of complex numbers such that the change in order of summation is justified, which is the case when $\{u_k\}$ terminates and when the double series converge absolutely. If $\{u_k\}$ is such that the sum over $k$ on the right side of (4.4) can be evaluated in terms of $q$-shifted factorials, then (4.4) yields a summation formula for a $q$-series. In particular, setting

$$u_k = \frac{(b,c,q^{-n};q)_k}{(q, aq, bcq^{-n}/a;q)_k} q^k$$

in (4.4), the sum over $k$ on the right side of (4.4) becomes a multiple of a terminating balanced $_3\phi_2$ series that can be summed by means of the $q$-Pfaff–Saalschütz summation formula (3.2), thus giving the sum of a terminating $_6W_5$ series

$$_6W_5(a;b,c,q^{-n};q,aq^{n+1}/bc) = \frac{(aq, aq/bc;q)_n}{(aq/b, aq/c;q)_n}, \qquad (4.5)$$

which reduces to (4.2) when $c = aq/b$. This formula was derived in BHS by first using (3.2) to derive the expansion formula (2.2.4) in BHS, and then applying (4.2) to a special case of the expansion formula. As in our derivation of (2.2) from (1.18), analytic continuation of (4.5) gives its nonterminating extension

$$_6W_5(a;b,c,d;q,aq/bcd) = \frac{(aq, aq/bc, aq/bd, aq/cd;q)_\infty}{(aq/b, aq/c, aq/d, aq/bcd;q)_\infty}, \qquad |aq/bcd| < 1. \quad (4.6)$$

This formula was obtained in §2.8 of BHS as a limit case of the summation formula (5.6) derived in §5.

In order to extend (4.6) to a summation formula for a $_6\psi_6$ series, it suffices to proceed as in the derivation of (2.8) from (2.2). Explicitly, replace the index of summation, call it $k$, in (4.6) by $k+n$, replace the parameters $a,b,c,d$ by $aq^{-2n}, bq^{-n}, cq^{-n}, dq^{-n}$, respectively, and manipulate the infinite products to

obtain that

$$_6\psi_6\left[\begin{array}{c} a/z, qa^{\frac{1}{2}}, -qa^{\frac{1}{2}}, b, c, d \\ zq, a^{\frac{1}{2}}, -a^{\frac{1}{2}}, aq/b, aq/c, aq/d \end{array}; q, \frac{azq}{bcd}\right]$$
$$= \frac{(aq, aq/bc, aq/bd, aq/cd, zq/b, zq/c, zq/d, q, q/a; q)_\infty}{(aq/b, aq/c, aq/d, q/b, q/c, q/d, zq, zq/a, azq/bcd; q)_\infty} \quad (4.7)$$

with $z = q^n$, $n = 0, 1, 2, \ldots$, where $0 < |azq/bcd| < 1$. Since the infinite series and the infinite product on the left and right sides, respectively, of (4.7) converge to analytic functions of $z$ when $|z| < \min(|1/q|, |a/q|, |bcd/aq|)$, by analytic continuation (4.7) holds when $|z| < \min(|1/q|, |a/q|, |bcd/aq|)$. Thus, setting $z = a/e$, we have derived Bailey's [1936] summation formula for a very-well-poised $_6\psi_6$ series

$$_6\psi_6\left[\begin{array}{c} qa^{\frac{1}{2}}, -qa^{\frac{1}{2}}, b, c, d, e \\ a^{\frac{1}{2}}, -a^{\frac{1}{2}}, aq/b, aq/c, aq/d, aq/e \end{array}; q, \frac{a^2q}{bcde}\right]$$
$$= \frac{(aq, aq/bc, aq/bd, aq/be, aq/cd, aq/ce, aq/de, q, q/a; q)_\infty}{(aq/b, aq/c, aq/d, aq/e, q/b, q/c, q/d, q/e, a^2q/bcde; q)_\infty}, \quad (4.8)$$

where $0 < |a^2q/bcde| < 1$. Also see the proofs referred to in §5.3 of BHS and, in particular, the extension of Ismail's [1977] proof of (2.8) to a proof of (4.8) presented in Askey and Ismail [1979].

Notice that, when $d = a^{\frac{1}{2}}$ and $e = -a^{\frac{1}{2}}$, (4.8) reduces to the following summation formula for a $_2\psi_2$ series that is different from the series considered in (2.10)

$$_2\psi_2(b, c; aq/b, aq/c; q, -aq/bc)$$
$$= \frac{(aq/bc; q)_\infty (aq^2/b^2, aq^2/c^2, q^2, aq, q/a; q^2)_\infty}{(aq/b, aq/c, q/b, q/c, -aq/bc; q)_\infty}, \qquad 0 < |aq/bc| < 1. \quad (4.9)$$

## 5. Additional summation and transformation formulas

In view of the important formulas that were derived in §4 from the expansion formula (4.4), which followed from the summation formula (4.2), it is of interest to investigate what formulas can be derived by replacing (4.2) by the more general summation formula (4.8). To introduce an integer parameter $k$, let us start by replacing $a, b, c, d, e$ in (4.8) by $aq^{2k}, bq^k, cq^k, dq^k, eq^{2k}$, respectively, to rewrite (4.8) in the form

$$_6\psi_6\left[\begin{array}{c} q^{k+1}a^{\frac{1}{2}}, -q^{k+1}a^{\frac{1}{2}}, bq^k, cq^k, dq^k, eq^{2k} \\ q^k a^{\frac{1}{2}}, -q^k a^{\frac{1}{2}}, aq^{k+1}/b, aq^{k+1}/c, aq^{k+1}/d, aq/e \end{array}; q, \frac{a^2 q^{1-k}}{bcde}\right]$$
$$= A \frac{(aq/b, aq/c, aq/d, be/a, ce/a, de/a; q)_k (a; q)_{2k}}{(b, c, d, bcde/a^2; q)_k (aq, e; q)_{2k}} \left(-\frac{bcd}{a}\right)^k q^{\binom{k}{2}} \quad (5.1)$$

with

$$A = \frac{(aq, aq/bc, aq/bd, aq/be, aq/cd, aq/ce, aq/de, q, q/a; q)_\infty}{(aq/b, aq/c, aq/d, aq/e, q/b, q/c, q/d, q/e, a^2q/bcde; q)_\infty}. \quad (5.2)$$

Using (5.1) with $A$ written as a (constant) function of $k$ and with $j$ as the index of summation in the $_6\psi_6$ series, we find that

$$A \sum_{k=-\infty}^{\infty} u_k = \sum_{k=-\infty}^{\infty} A\, u_k = \sum_{j=-\infty}^{\infty} \frac{(qa^{\frac{1}{2}}, -qa^{\frac{1}{2}}, b, c, d, e; q)_j}{(a^{\frac{1}{2}}, -a^{\frac{1}{2}}, aq/b, aq/c, aq/d, aq/e; q)_j} \left(\frac{a^2 q}{bcde}\right)^j$$
$$\times \sum_{k=-\infty}^{\infty} \frac{(bcde/a^2, eq^j, eq^{-j}/a; q)_k}{(be/a, ce/a, de/a; q)_k} u_k \qquad (5.3)$$

after replacing $j$ by $j - k$ and changing the order of summation, where it is assumed that $\{u_k\}_{k=-\infty}^{\infty}$ is a bilateral sequence of complex numbers such that the double series on the right side of (5.3) converges absolutely.

If we let $e = a$ in (5.3), which is equivalent to starting with the $_6\phi_5$ special case of (5.1), and set

$$u_k = \frac{(b, c, d, aq/ef; q)_k}{(q, aq/e, aq/f, bcd/a; q)_k} q^k, \qquad (5.4)$$

where at least one of numerator parameters $b, c, d, aq/ef$ is of the form $q^{-n}$, $n = 0, 1, 2\ldots$, so that the sequence $\{u_k\}$ has compact support, then the sum over $k$ on the right side of (5.3) becomes a terminating balance $_3\phi_2$ series that is summable via (3.2). This gives the transformation formula

$$_8W_7(a; b, c, d, e, f; q, a^2 q^2/bcdef)$$
$$= \frac{(aq, aq/bc, aq/bd, aq/cd; q)_\infty}{(aq/b, aq/c, aq/d, aq/bcd; q)_\infty} \,_4\phi_3\left[\begin{array}{c} b, c, d, aq/ef \\ aq/e, aq/f, bcd/a \end{array}; q, q\right] \qquad (5.5)$$

provided that the $_4\phi_3$ series on the right side terminates and $|a^2 q^2/bcdef| < 1$, so that the series on the left side converges. Watson [1929] used induction to prove (5.5) when both of the series terminate; and he commented that it should extend to the case when the series on the left converges and the series on the right terminates, which was subsequently proved by Bailey (see p. 70 in Bailey [1935]). Fifty years after Watson's paper was published, Askey and Ismail [1979] showed that (5.5) follows from Watson's terminating case by analytic continuation in the variable $z = 1/f$. Formula (5.5) can also be derived by observing that the proof of the expansion formula (2.5.2) in BHS extends to the case when $q^{-n}$ is replaced by a complex variable $f$ and $aq/bc$ is a negative integer power of $q$, so that using (4.6) to sum the resulting nonterminating very-well-poised $_6\phi_5$ series on the right side of the expansion formula yields (5.5). For a proof via orthogonal polynomials of (5.5) when both series terminate, see Andrews and Askey [1977].

When $b = q^{-n}$ and $a^2 q = bcdef$ the $_4\phi_3$ series in (5.5) reduces to a terminating balanced $_3\phi_2$ series that can be summed by (3.2) to derive, after replacing $f$ by $b$, Jackson's [1921] summation formula for a terminating $_8W_7$ series

$$_8W_7(a; b, c, d, e, q^{-n}; q, q) = \frac{(aq, aq/bc, aq/bd, aq/cd; q)_n}{(aq/b, aq/c, aq/d, aq/bcd; q)_n}, \qquad (5.6)$$

where $a^2 q^{n+1} = bcde$. Letting $n \to \infty$ in (5.6) gives (4.6). Note that the restrictions that the series in (5.6) and the series on the right side of (5.5) terminate cannot be removed by analytic continuation in a neighborhood of the origin of the $z$-plane by starting with the $z = q^n$, $n = 0, 1, 2, \ldots$, cases because there is no neighborhood of the origin in which all of the terms of the series are analytic. The $_4\phi_3$ series in (5.5) also reduces to a terminating balanced $_3\phi_2$ series that is summable by (3.2) when

$f = aq^{n+1}/e$ and $c = a/b$, which leads to a summation formula for a nonterminating $_8W_7$ series

$$_8W_7(a;b,a/b,d,e,aq^{n+1}/e;q,q^{1-n}/d)$$
$$= \frac{(q,aq,aq/bd,bq/d;q)_\infty}{(bq,aq/b,aq/d,q/d;q)_\infty} \frac{(aq/be,bq/e;q)_n}{(aq/e,q/e;q)_n}, \qquad |q^{1-n}/d| < 1, \qquad (5.7)$$

where $n = 0, 1, 2, \ldots$.

In order to use (5.7) to derive some new summation formulas, we first introduce a nonnegative integer parameter $k$ by replacing $a, b, c, d, e$ in (5.7) by $aq^{2k}, bq^k, cq^k, dq^k, eq^k$, respectively, to get

$$_8W_7(aq^{2k};bq^k,aq^k/b,dq^k,eq^k,aq^{k+n+1}/e;q,q^{1-k-n}/d)$$
$$= B_n \frac{(bq,aq/b,aq/d,aq/e,eq^{-n};q)_k}{(aq;q)_{2k}(d,e,aq^{n+1}/e;q)_k}(-d)^k q^{\binom{k}{2}+kn} \qquad (5.8)$$

with

$$B_n = \frac{(q,aq,aq/bd,bq/d;q)_\infty}{(bq,aq/b,aq/d,q/d;q)_\infty} \frac{(aq/be,bq/e;q)_n}{(aq/e,q/e;q)_n}. \qquad (5.9)$$

Then, proceeding as in the derivation of (5.3) with (5.1) replaced by (5.8), we obtain the expansion formula

$$B_n \sum_{k=0}^\infty u_k = \sum_{j=0}^\infty \frac{(a,qa^{\frac{1}{2}},-qa^{\frac{1}{2}},b,a/b,d,e,aq^{n+1}/e;q)_j}{(q,a^{\frac{1}{2}},-a^{\frac{1}{2}},aq/b,bq,aq/d,aq/e,eq^{-n};q)_j}\left(\frac{q^{1-n}}{d}\right)^j$$
$$\times \sum_{k=0}^j \frac{(aq^j,q^{-j};q)_k}{(b,a/b;q)_k} u_k, \qquad (5.10)$$

where $\{u_k\}_{k=0}^\infty$ is a sequence of complex numbers such that the double series on the right side of (5.10) converges absolutely.

If we let

$$u_k = \frac{(b,a/b,q^{-m};q)_k}{(q,aq/f,fq^{-m};q)_k} q^k \qquad (5.11)$$

in (5.10), where $m$ is a nonnegative integer, then the sums over $k$ on both the right and left sides of (5.10) are summable via (3.2), yielding the summation formula

$$_{10}W_9(a;b,a/b,d,e,aq^{n+1}/e,f,aq^{m+1}/f;q,q^{1-n-m}/d)$$
$$= \frac{(q,aq,aq/bd,bq/d;q)_\infty}{(bq,aq/b,aq/d,q/d;q)_\infty} \frac{(aq/be,bq/e;q)_n}{(aq/e,q/e;q)_n} \frac{(aq/bf,bq/f;q)_m}{(aq/f,q/f;q)_m}, \qquad (5.12)$$

where $n, m = 0, 1, 2, \ldots$, and $|q^{1-n-m}/d| < 1$. Moreover, by iterating this procedure, it follows that (5.7) extends to the family of summation formulas

$$_{6+2k}W_{5+2k}(a;b,a/b,d,e_1,\ldots,e_k,aq^{n_1+1}/e_1,\ldots,aq^{n_k+1}/e_k;q,q^{1-(n_1+\cdots+n_k)}/d)$$
$$= \frac{(q,aq,aq/bd,bq/d;q)_\infty}{(bq,aq/b,aq/d,q/d;q)_\infty}\prod_{j=1}^k \frac{(aq/be_j,bq/e_j;q)_{n_j}}{(aq/e_j,q/e_j;q)_{n_j}}, \qquad k=1,2,\ldots, \qquad (5.13)$$

where $n_1,\ldots,n_k$ are nonnegative integers and $|q^{1-(n_1+\cdots+n_k)}/d| < 1$.

For applications of (5.6) and additional summation, transformation, and expansion formulas, see BHS. The latest list of errata for BHS is available over the World Wide Web at: http://www.math.nwu.edu/preprints/gasper/index.html

## References

Andrews, G.E. [1969] *On a calculus of partition functions,* Pacific J. Math., **31**, 555–562.

Andrews, G.E. and Askey, R. [1977] *Enumeration of partitions: the role of Eulerian series and q-orthogonal polynomials,* Higher Combinatorics (M. Aigner, ed.), Reidel, Boston, Mass., pp. 3–26.

Andrews, G.E. and Askey, R. [1978] *A simple proof of Ramanujan's summation of the $_1\psi_1$,* Aequationes Math., **18**, 333–337.

Andrews, G. E. and Bressoud, D. M. [1984] *Identities in combinatorics III: Further aspects of ordered set sorting,* Discrete Math., **49**, 223-236.

Askey, R. [1980] *Ramanujan's extensions of the gamma and beta functions,* Amer. Math. Monthly, **87**, 346–359.

Askey, R. and Ismail, M.E.H. [1979] *A very-well-poised $_6\psi_6$,* Proc. Amer. Math. Soc., **77**, pp. 218–222.

Bailey, W.N. [1935] *Generalized Hypergeometric Series,* Cambridge University Press, Cambridge, reprinted by Stechert-Hafner, New York, 1964.

Bailey, W. N. [1936] *Series of hypergeometric type which are infinite in both directions,* Quart. J. Math. (Oxford), **7**, 105-115.

Bailey, W. N. [1941] *A note on certain q-identities,* Quart. J. Math. (Oxford), **12**, 173-175.

Berndt, B.C. [1993] *Ramanujan's theory of theta-functions,* Theta Functions From the Classical to the Modern (M. Ram Murty, ed.), CRM Proceedings & Lecture Notes, **1**, Amer. Math. Soc., Providence, R.I, pp. 1–63.

Bhatnagar, G. and Milne, S.C. [1995] *Generalized bibasic hypergeometric series and their $U(n)$ extensions,* to appear.

Cauchy, A.-L. [1843] *Mémoire sur les fonctions dont plusieurs valeurs sont liées entre elles par une équation linéaire, et sur diverses transformations de produits composés d'un nombre indéfini de facteurs,* C. R. Acad. Sci. Paris, T. XVII, p. 523, *Oeuvres de Cauchy,* 1$^{re}$ série, T. VIII, Gauthier-Villars, Paris, 1893, pp. 42–50.

Daum, J.A. [1942] *The basic analog of Kummer's theorem,* Bull. Amer. Math. Soc., **48**, 711–713.

Dougall, J. [1907] *On Vandermonde's theorem and some more general expansions,* Proc. Edin. Math. Soc., **25**, 114–132.

Gasper, G. [1989] *Summation, transformation, and expansion formulas for bibasic series,* Trans. Amer. Math. Soc., **312**, 257–277.

Gasper, G. and Rahman, M. [1990] *Basic Hypergeometric Series,* Encyclopedia of Mathematics and Its Applications, **35**, Cambridge University Press, Cambridge and New York.

Gasper, G. and Rahman, M. [1990a] *An indefinite bibasic summation formula and some quadratic, cubic, and quartic summation and transformation formulas,* Canad. J. Math., **42**, 1–27.

Gauss, C.F. [1813] *Disquisitiones generales circa seriem infinitam ...,* Comm. soc. reg. sci. Gött. rec., Vol. II; reprinted in *Werke* **3** (1876), pp. 123–162.

Goulden, I. [1985] *A bijective proof of the $q$-Saalschütz theorem,* Discrete Math., **57**, 39–44.

Hahn, W. [1949] *Über Polynome, die gleichzeitig zwei verschiedenen Orthogonalsystemen angehören,* Math. Nachr., **2**, 263-278.

Heine, E. [1847] *Untersuchungen über die Reihe ...,* J. reine angew. Math., **34**, 285–328.

Heine, E. [1878] *Handbuch der Kugelfunctionen, Theorie und Anwendungen,* Vol. 1, Reimer, Berlin.

Ismail, M.E.H. [1977] *A simple proof of Ramanujan's $_1\psi_1$ sum,* Proc. Amer. Math. Soc., **63**, 185–186.

Ismail, M.E.H. [1995] *The Askey–Wilson operator and summation Theorems,* Mathematical Analysis, Wavelets, and Signal Processing (M.E. Ismail, M.Z. Nashed, A.I. Zayed, and A.F. Ghaleb, eds.), Contempory Mathematics, **190**, Amer. Math. Soc., Providence, R.I, pp. 171–178.

Jackson, F.H. [1910] *Transformations of $q$-series,* Messenger of Math., **39**, 145–153.

Jackson, F.H. [1921] *Summation of $q$-hypergeometric series,* Messenger of Math., **50**, 101–112.

Jacobi, C.G.J. [1829] *Fundamenta Nova Theoriae Functionum Ellipticarum,* Regiomonti. Sumptibus fratrum Borntraeger; reprinted in Gesammelte Werke, **1** (1881), 49–239, Reimer, Berlin.

Jacobi, C. G. J. [1846] *Über einige der Binomialreihe analoge Reihen,* J. reine angew. Math., **32**, 197-204; reprinted in Gesammelte Werke, **6** (1881), 163-173, Reimer, Berlin.

Pfaff, J.F. [1797] *Observationes analyticae ad L. Euler Institutiones Calculi Integralis,* Vol. IV, Supplem. II et IV, Historia de 1793, Nova acta acad. sci. Petropolitanae, **11** (1797), pp. 38–57.

Rahman, M. [1990] *Some extensions of the beta integral and the hypergeometric function,* Orthogonal Polynomials (P. Nevai, ed.), Kluwer Academic Publishers, pp. 1–27.

Saalschütz, L. [1890] *Eine Summationsformel,* Zeitschr. Math. Phys., **35**, 186–188.

Slater, L.J. [1966] *Generalized Hypergeometric Functions,* Cambridge University Press, Cambridge.

Watson, G.N. [1929] *A new proof of the Rogers-Ramanujan identities,* J. London Math. Soc., **4**, 4–9.

Wilf, H.S and Zeilberger, D. [1992] An algorithmic proof theory for hypergeometric (ordinary and "$q$") multisum/integral identities, **108**, 575–633.

Zeilberger, D. [1987] A $q$-Foata proof of the $q$-Saalschütz identity, Europ. J. Combinatorics, **8**, 461–463.

E-mail address:   george@math.nwu.edu